\documentclass[11pt,leqno]{article}  
\usepackage{amsmath,amssymb}  
\usepackage{epsfig}  

\begin{document}  

\newtheorem{defi}{Definition}[section]  
\newtheorem{rema}{Remark}[section]  
\newtheorem{prop}{Proposition}[section]  
\newtheorem{lem}{Lemma}[section]  
\newtheorem{theo}{Theorem}[section]  
\newtheorem{cor}{Corollary}[section]  
\newtheorem{conc}{Conclusion}[section]  
  
\author{P. Germain}  
  
\title{Space-time resonances}

\maketitle
  
\begin{abstract} 
This article is a short exposition of the space-time resonances method. It was introduced by Masmoudi, Shatah, and the author, in order to understand global existence for nonlinear dispersive equations, set in the whole space, and with small data. The idea is to combine the classical concept of resonances, with the feature of dispersive equations: wave packets propagate at a group velocity which depends on their frequency localization. The analytical method which follows from this idea turns out to be a very general tool.
\end{abstract}

This expository article aims at presenting the space-time resonance method, as it is developped by Masmoudi, Shatah, and the author~\cite{GMS1}~\cite{GMS2}~\cite{GMS3}, and the author~\cite{G}.

\section{Introduction: the problem of global existence for small data}  

\label{cardinal}

Given the Cauchy problem
\begin{itemize}
\item for a nonlinear dispersive equation
\item set in the whole space $\mathbb{R}^d$ ($d\geq 1$)
\item with Cauchy data which are \underline{small, smooth and localized},
\end{itemize}
we want to answer the following question: does there exist a global in time solution? Does it scatter?

\bigskip

This question can be anwered easily if the nonlinearity corresponds to a sufficiently high power, depending on the decay given by the linear part of the equation. ``Sufficiently high'' can often be quantified by the Strauss exponent
$$
\gamma(d) \overset{def}{=} \frac{1}{2} + \frac{1}{d} + \sqrt{\left(\frac{1}{2} + \frac{1}{d}\right)^2 + \frac{2}{d}}.
$$
Assume the nonlinearity is of the form $f(u)$, for a function $f$ such that $|f'(u)|\lesssim |u|^{p-1}$. If the linear part is the Schr\"odinger equation, then global existence for small, smooth and localized data follows for $p > \gamma(d)$, see Strauss~\cite{S}; if it is the wave equation, one needs to ask $p>\gamma(d-1)$ (due to the weaker decay), see John~\cite{J} and Georgiev, Lindblad and Sogge~\cite{GLS}.

\bigskip

Thus, caricaturing a little, global existence follows from dispersive estimates, regardless of the precise structure of the nonlinear term, provided its homogeneity is above the Strauss exponent. Below it, it can be shown that the structure of the nonlinearity matters, since resonances start playing a role. To this we now turn.

\section{The notion of space-time resonance}

\label{swan}

\subsection{A simple case to begin with}

\label{bluejay}

To make the exposition as simple as possible, we consider for most of this paper the case of the scalar equation
\begin{equation}
\label{baldeagle}
\left\{ \begin{array}{l}
i \partial_t u + P(D) u = N(u) \\
u(t=0) = u_0,
\end{array} \right.
\end{equation}
where $u(t,x)$ is a complex-valued function of the time variable $t \in \mathbb{R}$ and of the space variable $x \in \mathbb{R}^d$, $P(D)$ is a real Fourier multiplier with symbol $P(\xi)$, the dispersion relation, and
$$
N = N_{\epsilon_1 \epsilon_2}, \;\;\; \mbox{with} \; \epsilon_1,\epsilon_2 = + \;\mbox{or}\;-
$$
and
$$
N_{++}(u) = u^2 \;\;\;\;,\;\;\;\; N_{--}(u) = \bar u^2\;\;\;\;\mbox{and}\;\;\;\;N_{+-}(u)=u\bar u.
$$
As underlined above, we always assume that $u_0$ is small, smooth and localized.

\subsection{The viewpoint of stationary phase}

The idea is, after rewriting the equation in a convenient way, to analyze the nonlinear term by the stationary phase method; thus we want to isolate all the oscillations in a unique factor. Since the solution $u$ displays oscillations itself, the first step is to take as the new unknown function
$$
f(t) \overset{def}{=} e^{-itP(D)} u(t).
$$
Write then Duhamel's formula for the Fourier transform $\widehat{f}$ of $f$: it gives, as follows by a small computation,
\begin{equation}
\label{eagle}
\widehat{f}(t,\xi) = \widehat{u_0}(\xi) + \int_0^t \int e^{is\phi(\xi,\eta)} \widehat{f_{\epsilon_1}}(s,\eta) \widehat{f_{\epsilon_2}}(s,\xi-\eta)\,d\eta\,ds,
\end{equation}
where we denote $f_+\overset{def}{=}f$, $f_-\overset{def}{=}\bar f$, and furthermore
\begin{equation*}
\phi(\xi,\eta) \overset{def}{=} P(\xi)- \epsilon_1 P(\eta)- \epsilon_2 P(\xi-\eta) 
\end{equation*}
In order to prove global existence and scattering, it is desirable to control the solution uniformly in time; but the domain of the integral in~(\ref{eagle}) over $s$ grows with $t$. This can be compensated by the oscillating factor $e^{is\phi(\xi,\eta)}$, excepts on the sets where the phase is stationary in either of the integration variables:
\begin{itemize}
\item On the set of \underline{time resonances} $\displaystyle \mathcal{T} \overset{def}{=} \{(\xi,\eta) \;\;\mbox{such that}\;\;\phi(\xi,\eta) = 0$\}, the phase is stationary in $s$.
\item On the set of \underline{space resonances} $\displaystyle \mathcal{S} \overset{def}{=} \{(\xi,\eta) \;\;\mbox{such that}\;\;\partial_\eta \phi(\xi,\eta) =0$\}, the phase is stationary in $\eta$.
\item On the set of \underline{space-time resonances} $\displaystyle \mathcal{R} \overset{def}{=} \mathcal{S} \cap \mathcal{T}$, the phase is stationary in $s$ and $\eta$.
\end{itemize}
The key idea is that the sets $\mathcal{T}$, $\mathcal{S}$, and even more $\mathcal{R}$, are the obstructions to a linear behaviour of $u$ for large time.

\subsection{The physical meaning}

It is striking that the elementary stationary phase analysis above captures deep physical properties of the equation. In order to demonstrate this, notice that, in first approximation, the nonlinear term in~(\ref{baldeagle}) will simply be the product of solutions of the linear equation
\begin{equation}
\label{duck}
\epsilon_j i \partial_t u + P(D) u = 0\;\;\;\;\mbox{with $j=1,2$}.
\end{equation}

\subsubsection{Time resonances}

To understand the physical meaning of time resonances, consider plane wave solutions of~(\ref{duck}): they read for $j=1,2$ 
$$
u^j_{pw}(t,x) = e^{i\left( x\cdot\xi - \epsilon_j tP(\xi_j) \right)}
$$
(for some given frequencies $\xi_1$, $\xi_2$).
The product of $u_{pw}^1$ and $u_{pw}^2$ is
$$
e^{i\left( x\cdot\left[\xi_1 + \xi_2\right] - t [\epsilon_1 P(\xi_1) + \epsilon_2 P(\xi_2)] \right)}
$$
If the expression above is a plane wave, then the interaction is resonant in the classical sense (as for ODEs or dynamical systems). But this is the case if and only if $\phi (\xi_1 + \xi_2,\xi_2)=0$!

\bigskip

To summarize: the plane waves $u^1_{pw}$ and $u^2_{pw}$ are resonant (in the classical sense) if and only if $(\xi_1+\xi_2,\xi_2) \in \mathcal{T}$.

\subsubsection{Space resonances}

Space resonances can be understood by considering the interaction of wave packets, which are the localized version of plane waves. Consider two wave packets $u^1$ and $u^2$ solving~(\ref{duck}) with given data:
$$
\left\{ \begin{array}{l} \epsilon_j i \partial_t u^j + P(D) u^j = 0 \\ u(t=0)=u_0^j. \end{array} \right.
$$
We next pick data $u_0^j$ which are localized around $0$ in space and $\xi_j$ in frequency. It is well known that the $u^j$ will propagate at the group velocity given by $- \epsilon_j \nabla P (\xi_j)$; in other words, $u^j$ will be (essentially) localized around the set $x \sim - \epsilon_j\nabla P (\xi_j) t $. It is clear that these two wave packets can, for large time, have a non-negligible interaction only if their group velocities agree, meaning $\epsilon_1 \nabla P(\xi_1) = \epsilon_2 \nabla P(\xi_2)$. But this is eqivalent to $\partial_\eta \phi(\xi_1 + \xi_2,\xi_2)=0$!

\bigskip

To summarize: two wave packets with frequencies $\xi_1$ and $\xi_2$ share the same space-time support if and only if $(\xi_1+ \xi_2,\xi_2) \in \mathcal{S}$.

\section{The analytical method}

\label{robin}

Now that we achieved a good intuitive understanding of the problem, we turn to the analytical work: deriving estimates.

\subsection{Multilinear cutoff in frequency space}

\label{sparrow}

In order to get estimates, the idea will be to distinguish regions in the $(\xi,\eta)$ plane depending on the oscillatory properties of $e^{is\phi}$: for instance, away from $\mathcal{S}$ one can take advantage of oscillations in $s$, and away from $\mathcal{T}$ of oscillations in $\eta$.

Thus one needs to cut-off in the $(\xi,\eta)$ plane; this naturally leads to the pseudo-product operators of Coifman and Meyer~\cite{CM}. Given a symbol $m(\xi,\eta)$, the pseudo-product with this symbol reads
$$
B_m (f,g) \overset{def}{=} \mathcal{F}^{-1} \int m(\xi,\eta) \widehat{f}(\eta) \widehat{g}(\xi-\eta) \,d\eta.
$$
Choosing $m$ equal to a cut-off function enables one to consider separately different regions of $(\xi,\eta)$ space.

\label{parrot}

\subsection{Away from $\mathcal{T}$: normal form}

Suppose that a cut-off as in the previous paragraph enabled us to focus on a region where $\phi$ does not vanish, and let us forget about the cut-off function in order to simplify the expressions. The most natural thing to do in order to take advantage of the oscillations is to integrate by parts using the identity
$$
\frac{1}{i\phi} \partial_s e^{is \phi} = e^{is\phi}.
$$
This turns~(\ref{eagle}) into
$$
\widehat{f}(t,\xi) \overset{def}{=} \widehat{u_0}(\xi) + \int \frac{1}{i\phi} e^{it\phi(\xi,\eta)} \widehat{f}(t,\eta) \widehat{f}(t,\xi-\eta)\,d\eta - \int_0^t \int 
\frac{1}{i \phi} e^{is\phi(\xi,\eta)} \partial_s \widehat{f}(s,\eta) \widehat{f}(s,\xi-\eta)\,d\eta\,ds
$$
plus other symmetric or less important terms, which we do not detail. This transformation is favorable as far as estimates are concerned: the second term in the above right-hand side is independent of $t$, whereas the third is actually cubic in $u$ (or $f$), due to the relation $\partial_t f = e^{-itP(D)} N(u)$.

\bigskip

As can easily be seen, this procedure is nothing but a reformulation of the normal form method, first used in a dispersive PDE context by Shatah~\cite{S}.

\subsection{Away from $\mathcal{S}$: vector fields}

Just as above, assume that a cut-off, which we then forget, enabled us to focus on a region where $\partial_\eta \phi$ does not vanish. Use then the identity
$$
\frac{1}{is |\partial_\eta \phi|^2} \partial_\eta \phi \cdot \partial_\eta e^{is \phi} = e^{is\phi}
$$
in order to integrate by parts in~(\ref{eagle}), yielding 
$$
\widehat{f}(t,\xi) \overset{def}{=} \widehat{u_0}(\xi) - \int_0^t 
\int e^{is\phi(\xi,\eta)} \frac{1}{is|\partial_\eta \phi|^2} \partial_\eta \phi \cdot \partial_\eta \widehat{f}(s,\eta) \widehat{f}(s,\xi-\eta)\,d\eta\,ds
$$
plus other symmetric or easier terms. What is gained by this transformation is clear: an extra factor $\frac{1}{s}$ in the bilinear term which helps controlling it for large $s$.

\bigskip

In order to recast the above manipulation in a more familiar setting, we need to come back to the original unknown function $u$. It turns out that
$$
e^{itP(D)} \mathcal{F}^{-1} \partial_\eta \widehat{f}(t) = J u(t),
$$
where $J = ix - t \left[\nabla P\right] (D)$. This is precisely the operator used in the vector fields method for this type of equation! The idea of the vector fields method, first introduced by Klainerman~\cite{K} in the context of the nonlinear wave equation is the following: do not estimate only $u$ in various Sobolev-type norms, but also $Lu$, where $L$ belongs to a family of vector fields which commute with the linear part of the equation.

Thus our integration by parts essentially corresponds to a new point of view on the vector fields method.

\subsection{What about $\mathcal{R}$?}

We saw in the previous two paragraphs that the course of action is fairly clear if one is away from $\mathcal{R}$, and we want to discuss now how a neighbourhood of $\mathcal{R}$ can be treated. It seems that a fairly good idea is to add a cut-off function
$$
\chi\left(t^\delta \operatorname{dist}((\xi,\eta),\mathcal{R})\right),
$$
where $\chi \in \mathcal{C}^\infty_0$ is equal to $1$ in a neighbourhood of zero, and $\delta>0$; this amounts to cutting off to a distance $t^{-\delta}$ of $\mathcal{R}$. This can be exploited as follows:
\begin{itemize}
\item For frequencies within $t^{-\delta}$ of $\mathcal{R}$, one can take advantage of the shrinking size of this set (it is reasonable to assume that $\mathcal{R}$ has measure zero).
\item For frequencies lying more than $t^{-\delta}$ away from $\mathcal{R}$, one can integrate by parts, either in $s$, or in $\eta$, and hope to balance the loss due to $\phi$ and $\partial_\eta \phi$ being closer to $0$ as $\mathcal{R}$ is approached.
\end{itemize}

\subsection{Conclusion: comparison with classical methods}

As we just saw, partitionning the $(\xi,\eta)$ space allows one to use (equivalent formulations of) either the normal form, or the vector fields method, depending on the elementary wave interaction which is being considered. This is often a decisive improvement.

Furthermore, the approach that was just sketched is very systematic and adaptable; it does not rely on delicate algebraic identities as is sometimes the case for the vector fields method.

However, in the case of the wave equation, it seems that the vector fields method (which was originally introduced in that context by Klainerman~\cite{K}) provides a more direct proof of the classical theorems. It might be linked to the equation being second order in time.

\section{Some other aspects of space-time resonances}

\subsection{Outcome, source and separation}

When dealing with resonances, an important question is: can resonant interactions form a loop? To be more precise, define the projections of $\mathcal{R}$ on $\xi$, $\eta$, and $\xi-
\eta$:
\begin{equation*}
\begin{split}
& \mbox{outcome frequencies} \overset{def}{=} \{ \xi\;\mbox{such that}\;(\xi,\eta) \in\mathcal{R} \} \\
& \mbox{source frequencies} \overset{def}{=} \{ \eta\;\mbox{such that}\;(\xi,\eta)\;\mbox{or}\;(\xi,\xi-\eta) \in \mathcal{R} \}.
\end{split}
\end{equation*}
The terminology is clear: outcome frequencies can be seen as being produced by a space-time resonance, whereas source frequencies can be at the origin of one. If
$$
\mbox{outcome frequencies}\;\neq\;\mbox{source frequencies},
$$
resonances are said to be separated. It is clear that this condition will help in controlling the equation: indeed, if it holds true, the worst interactions (the resonant ones) cannot be fed by the frequencies which are least controlled (outcome frequencies).

See~\cite{G} for an instance of this phenomenon, and of its utilization.

\subsection{Cancellations}

\subsubsection{Null form}

A null form is a structure of the nonlinear term such that resonant interactions are canceled; in other words, the nonlinearity vanishes on the points where the phase is stationary.

To be more precise, assume that the nonlinear term $N(u)$ is given by a pseudo-product: $N(u) = B_m(u,u)$. If one tries to perform directly (without localizing) the manipulations of Section~\ref{robin}, there appears pseudo-products with the symbols $\displaystyle \frac{m(\xi,\eta)}{\phi(\xi,\eta)}$ and $\displaystyle \frac{m(\xi,\eta) \partial_\eta\phi}{|\partial_\eta \phi|^2}$.

Are these new pseudo-products bounded? This of course what matters, more than the vanishing of $\phi$ or $\partial_\eta \phi$. Caricaturing a little, if $m$ vanishes to a higher order than $\phi$ or $\partial_\eta \phi$, then the nonlinear term can be controlled.

\subsubsection{Vanishing of $\partial_\xi \phi$}

The most important features of $\phi$ are its zero set, $\mathcal{T}$, and that of $\partial_\eta \phi$, $\mathcal{S}$. But the zero set of $\partial_\xi \phi$ can also turn out to be useful. More specifically: if the zero set of $\partial_\xi \phi$ contains $\mathcal{T}$, $\mathcal{S}$, or $\mathcal{R}$, estimates can be derived more easily. This property is for instance used in~\cite{GMS2} and~\cite{GMS3}.

\bigskip

The analytical explanation is the following: suppose you want to estimate $xf$; it corresponds in Fourier space to $\partial_\xi \widehat{f}$. 
The worst term one obtains when applying $\partial_\xi$ to~(\ref{eagle}) contains an additional $s$ factor: it is
$$
\int_0^t \int is\partial_\xi \phi(\xi,\eta) e^{is\phi(\xi,\eta)} \widehat{f}(s,\eta) \widehat{f}(s,\xi-\eta)\,d\eta\,ds.
$$
If now $\partial_\xi \phi$ vanishes only on $\mathcal{T}$, $\mathcal{S}$, or $\mathcal{R}$, this is very reminiscent of a null form structure!

\subsection{The general case} 

\label{pinguin}

\subsubsection{General setting}

We want to show in this paragraph how the Cauchy problem for a general nonlinear dispersive equation can be dealt with very similarly to the simple case of Subsection~\ref{bluejay}. Thus consider the Cauchy problem
\begin{equation*}
\left\{ \begin{array}{l}
i \partial_t u + P(D) u = N(u) \\
u(t=0) = u_0
\end{array} \right.
\end{equation*}
where $u(t,x)$ is an $\mathbb{R}^n$-valued function of time $t \in \mathbb{R}$ and space $x \in \mathbb{R}^d$; $N$ is a smooth nonlinear function of $u$, $\bar u$; and we can assume without loss of generality that $P$ is given by a diagonal matrix of real Fourier multipliers
$$
P(D) = \operatorname{diag}(P_1(D),P_2(D)\dots P_n(D)).
$$
As we now show, the oscillatory phase approach which has been outlined above can be adapted in a straightforward fashion to this case.

\subsubsection{Treating the more general nonlinearity}

First, taking for $N$ polynomials of order higher than 2 only changes a few details to the manipulations performed in Section~\ref{swan}; the generalization to $N$ being a pseudo-product operator (see paragraph~\ref{parrot}) is also immediate. This suffices to treat the general case since a nonlinearity which is smooth and translation invariant can be expanded as a series of multilinear pseudo-product operators. The first summands of these series can be subjected to the space-time resonance approach; higher order summands correspond to a large power, and thus can be treated easily, as noticed in Section~\ref{cardinal}.

This method for treating a general nonlinearity is used in~\cite{GMS3}.

\subsubsection{Dealing with the vectorial structure}

Following the path laid out in Section~\ref{bluejay}, phases appear which correspond to interactions between the different scalar waves: for instance, phases corresponding to quadratic interactions look like
$$
\phi (\xi,\eta) = P_i(\xi) \pm P_j(\eta) \pm P_k(\xi-\eta)\;\;\;\;\mbox{where} \;i,j,k \in \{1,2,\dots,n\}.
$$
Aside from this modification, all the above discussions remains valid.

\subsection{A few examples of resonant sets}

It is instructive to examine what resonant sets look like for various equations; we still consider only quadratic interactions.

\subsubsection{Scalar equation with homogeneous dispersion relation}

The setting is that of Section~\ref{bluejay}, with the additional assumption that $P(\xi) = |\xi|^\alpha$, for some $\alpha>0$.
\begin{itemize}
\item If $0<\alpha<1$, $\mathcal{S}$ is always a linear subspace; in the cases $++$ and $--$, $\mathcal{T}=\{(0,0)\}$; in the case $+-$, $\mathcal{R} = \{\xi=0\}$.
\item If $\alpha = 1$, $\mathcal{S}$, $\mathcal{R}$ and $\mathcal{T}$, when non-empty, are of the form $\{(\xi,\eta),\xi\;\mbox{and}\;\eta\;\mbox{positively colinear}\}$.
\item If $\alpha>1$, $\mathcal{S}$ is always a linear subspace; in the case $++$, $\mathcal{T}=\{(0,0)\}$; in the case $--$, $\mathcal{R}=\{(0,0)\}$; and in the case $+-$ $\mathcal{R} = \{\xi=0\}$.
\end{itemize}
To summarize: the case $0<\alpha<1$ is the best behaved; the case $\alpha=1$ (wave equation) is very particular; and for $\alpha \neq 1$, linear subspaces tend to play a key role.

\subsubsection{General case}

We adopt the general setting of Subsection~\ref{pinguin}, and make the additional assumption that the dispersion relations $P_i(\xi)$ only depend on $|\xi|$.
Then, as explained in Subsection~\ref{pinguin}, quadratic phases have the form
$$
\phi(\xi,\eta) = P_1(|\xi|)\pm P_2(|\eta|) \pm P_3(|\xi-\eta|)
$$
for some functions $P_1$, $P_2$ and $P_3$.
A very simple computation gives then, under generic assumptions, that $\mathcal{R}$ is either empty, or of the form $\{|\xi|=R\;\mbox{and}\; \eta = \lambda \xi \}$, for real constants $R$ and $\lambda$.

\section{Related questions}

\subsection{Local well-posedness}

It is natural that the two building blocks of our method: time resonances (linked to classical resonances) and space resonances (linked to group velocity) also play an important role for local well-posedness problems. They appear however under a slightly different disguise, which we now describe.

\bigskip

A very useful tool is provided by $X^{s,b}$ spaces; they  have been introduced by Bourgain~\cite{Bo} and Klainerman and Machedon~\cite{KM}, in order to investigate local well-posedness of nonlinear dispersive equations for low regularity data. They are given by the norm
$$
\| u \|_{X^{s,b}} \overset{def}{=} \int \int |\xi^{2s}| |\tau-P(\xi)|^{2b} \left| \widetilde{\chi u}(\tau,\xi) \right|^2 \,d\xi \,d\tau
$$
where $\chi$ is a cut-off function localizing in time to a neighbourhood of zero ($u$ is controlled only for small times); $P$ is the dispersion relation of the equation under study; and $\widetilde{\cdot}$ denotes the space-time Fourier transform.

These spaces can be understood as measuring how close to a solution of the linear equation $u$ is, for small times. For large $b$ at least, they therefore reduce the nonlinear problem to understanding the interaction of linear waves; it is then natural that resonances will play a key role. A very nice and more detailed explanation can be found in~\cite{T}.

\bigskip

Another instance where spaces resonances come up is the following: Strichartz inequalities can be improved by considering multilinear estimates, and taking the argument functions to have different frequency supports. We borrow an instance of this phenomenon from Bourgain~\cite{Bo1}: if $\psi_1$ is localized in frequency around $M_1$, and $\psi_2$ around $M_2$, then
$$
\left\|e^{it\Delta} \psi_1 e^{it\Delta} \psi_2 \right\|_{L^2(\mathbb{R} \times \mathbb{R}^2)} \lesssim \left( \frac{M_1}{M_2} \right)^{1/2} \left\|\psi_1\right\|_{L^2(\mathbb{R}^2)} \left\|\psi_2\right\|_{L^2(\mathbb{R}^2)},
$$
which refines the Strichartz inequality 
$$
\left\|e^{it\Delta} \psi \right\|_{L^2(\mathbb{R} \times \mathbb{R}^2)} \lesssim \|\psi\|_{L^2(\mathbb{R}^2)}.
$$
This improvement, very useful to study nonlinear interactions, of course follows from the difference in group velocities of wave packets with different frequency supports.

\subsection{Multilinear oscillatory operators}

Operators of the type
\begin{equation}
\label{seagull}
(f,g) \mapsto \mathcal{F}^{-1} \int_0^t \int m(\xi,\eta) e^{is\phi(\xi,\eta)} \widehat{f}(\xi-\eta) \widehat{g}(\eta) \,d\eta\,ds
\end{equation}
(or even the $n$-linear version of these) obviously played a key role in the above discussion. Curiously enough, this class of operators has not received much attention as such.
It would be interesting to understand its boundedness properties, between say Lebesgue or weighted Lebesgue spaces, with a bound which may depend on $t$.

\bigskip

A moment of reflection shows that this question is related to the boundedness of pseudo-product operators (defined in Section~\ref{sparrow}) with singular (say, discontinuous along a manifold) symbols. Boundedness properties of such operators are far from being understood; a notable exception is the bilinear Hilbert transform, see the celebrated article of Lacey and Thiele~\cite{LT}.

\bigskip

Bernicot and Germain~\cite{BG} studied the simplified case where the time integral is removed in~(\ref{seagull}):
$$
(f,g) \mapsto \mathcal{F}^{-1} \int m(\xi,\eta) e^{is\phi(\xi,\eta)} \widehat{f}(\xi-\eta) \widehat{g}(\eta) \,d\eta
$$
(which amounts to ignoring time resonances). They were able to derive optimal bounds for generic $\phi$; notice however that the problem becomes more simple if $\phi$ splits as a sum of functions of $\xi$, $\eta$, and $\xi-\eta$, which is the case for PDE applications.

\bigskip  
  
\bigskip  
  
\bigskip  
  
Pierre GERMAIN  
  
{\sc Courant Institute of Mathematical Sciences  
  
New York University   
  
New York, NY 10012-1185  
  
USA}  
  
\bigskip  
    
{\tt pgermain@math.nyu.edu}  
  
\end{document}